\RequirePackage{fix-cm}
\documentclass[smallextended]{svjour3}       % onecolumn (second format)
\smartqed  % flush right qed marks, e.g. at end of proof
\usepackage{graphicx}
\usepackage{amsmath}
\usepackage{amsfonts}
\usepackage{cite}
\journalname{Results in Mathematics}
\begin{document}

\title{A note on some reduction formulas for the incomplete beta function
and the Lerch transcedent%\thanks{Grants or other notes
}

\titlerunning{Reduction formulas for incomplete beta and Lerch transcendent}        % if too long for running head

\author{J.L. Gonz\'{a}lez-Santander}


\institute{J.L. Gonz\'{a}lez-Santander \at
              Department of Mathematics. Universidad de Oviedo. \\
              Federico Garc\'{i}a Lorca 18, 33007 Oviedo, Spain. \\
              Tel.: +34-985-10-3338\\
              \email{gonzalezmarjuan@uniovi.es}           %  \\
}

\date{Received: date / Accepted: date}

\maketitle

\begin{abstract}
We derive new reduction formulas for the incomplete beta function and the
Lerch transcendent in terms of elementary functions. As an application, we
calculate some new integrals. Also, we use these reduction formulas to test
the performance of the algorithms devoted to the numerical evaluation of the
incomplete beta function.
\end{abstract}

\section{Introduction}

The origins of the beta function $\mathrm{B}\left( \nu ,\mu \right) $ go
back to Wallis' attempt of the calculation of $\pi $ \cite{Wallis}. For this
purpose, he evaluated the integral
\begin{equation}
\mathrm{B}\left( \nu ,\mu \right) =\int_{0}^{1}t^{\nu -1}\left( 1-t\right)
^{\mu -1}dt,  \label{beta_wallis}
\end{equation}%
where $\nu $ and $\mu $ are integers or $\mu =1$ and $\nu $ is rational.
Moreover, Wallis suggested that \cite[p. 4]{Andrews}
\begin{equation*}
\frac{\pi }{4}=\frac{1}{2}\int_{0}^{1}t^{-1/2}\left( 1-t\right) ^{1/2}dx=%
\frac{1}{4}\lim_{n\rightarrow \infty }\left( \frac{2\cdot 4\cdot 6\ \cdots \
2n}{1\cdot 3\cdot 5\ \cdots \ \left( 2n-1\right) }\frac{1}{\sqrt{n}}\right)
^{2}.
\end{equation*}

This result may have led Euler to consider the integral (\ref{beta_wallis})\
for $\nu $ and $\mu $ not neccesarily integers and its relation to the gamma
function. In fact, Euler derived the following relation between the beta and
gamma functions \cite[Eqn. 1.1.13]{Andrews}:
\begin{equation*}
\mathrm{B}\left( \nu ,\mu \right) =\frac{\mathrm{\Gamma} \left( \nu \right) \mathrm{\Gamma}
\left( \mu \right) }{ \mathrm{\Gamma} \left( \nu +\mu \right) }.
\end{equation*}

A natural generalization of the beta function is the incomplete beta
function, defined as \cite[Eqn. 8.17.1]{DLMF}%
\begin{equation*}
\mathrm{B}\left( \nu ,\mu ,z\right) =\int_{0}^{z}t^{\nu -1}\left( 1-t\right)
^{\mu -1}dt,\qquad 0\leq z\leq 1,\ \nu ,\mu >0,
\end{equation*}%
where it is straightforward to continue analytically to complex values of $%
\nu $, $\mu $, and $z$.

Many applications have been developed over time regarding the $\mathrm{B}%
\left( \nu ,\mu ,z\right) $ function. For instance, in statistics it is used
extensively as the probability integral of the beta distribution \cite[p.
210-275]{Johnson}. Also, it appears in statistical mechanics for Monte Carlo
sampling \cite{Kofke}, in the analysis of packings of granular objects \cite%
{Prellberg}, and in growth formulas in cosmology \cite{Hamilton}. Therefore,
in order to evaluate the $\mathrm{B}\left( \nu ,\mu ,z\right) $ function, it
is quite interesting to have reduction formulas to simplify its computation,
both symbolically and numerically. For instance, when $\mu =m+1$\ is a
positive integer (i.e. $m=0,1,2,\ldots $), we have the following reduction
formula in terms of elementary functions \cite[Eqn. 58:4:3]{Atlas}
\begin{equation}
\mathrm{B}\left( \nu ,m+1,z\right) =z^{\nu }\sum_{k=0}^{m}\left(
\begin{array}{c}
m \\
k%
\end{array}%
\right) \frac{\left( -z\right) ^{k}}{k+\nu }.  \label{Reduction_Atlas}
\end{equation}

However, when $\mu =0$, the incomplete beta function is given in terms of
the Lerch transcendent \cite[Eqn. 58:4:4]{Atlas}%
\begin{equation}
\mathrm{B}\left( \nu ,0,z\right) =z^{\nu }\mathrm{\Phi} \left( z,1,\nu \right) ,\quad
\nu >0,  \label{B=Phi}
\end{equation}%
where the Lerch transcendent is defined as \cite[Eqn. 1.11(1)]{Erdelyi}
\begin{equation}
\mathrm{\Phi} \left( z,s,\nu \right) =\sum_{k=0}^{\infty }\frac{z^{k}}{\left( k+\nu
\right) ^{s}},\quad \left\vert z\right\vert <1,\ \nu \neq 0,-1,-2,\ldots
\label{Lerch_def}
\end{equation}

It is worth noting that (\ref{Reduction_Atlas})\ can be proved by induction
from (\ref{B=Phi}) and (\ref{Lerch_def}), applying the connection formula
\cite[Eqn. 58:5:3]{Atlas}:
\begin{equation*}
\mathrm{B}\left( \nu ,\mu ,z\right) =\mathrm{B}\left( \nu +1,\mu ,z\right) +%
\mathrm{B}\left( \nu ,\mu +1,z\right) .
\end{equation*}

Nevertheless, reduction formulas for $\mathrm{B}\left( \nu ,0,z\right) $\
when $\nu $ is a rational number do not seem to be reported in the most
common literature. The aim of this note is just to provide such reduction
formulas in terms of elementary functions. As an application, we will
calculate some new integrals in terms of elementary functions. Also, we will
check that the numerical evaluation of the incomplete beta function is
improved with these reduction formulas.

This paper is organized as follows. Section \ref{Section: Main Results}
derives reduction formulas for $\mathrm{B}\left( \nu ,0,z\right) $, both for
$\nu $ positive rational, as well as negative rational. Particular cases for
$\nu $ non-negative integer or $\nu $ half-integer are also derived. In
Section \ref{Section: Applications}, we will apply the reduction formulas
derived in Section \ref{Section: Main Results} to calculate some integrals
which do not seem to be reported in the most common literature. Further, for
particular values of the parameters, the symbolic computation of these
integrals is quite accelerated by using the aforementioned reduction
formulas. Also, we will use these reduction formulas to numerically test the
performance of the algorithm provided in MATHEMATICA to
compute the incomplete beta function.

\section{Main results\label{Section: Main Results}}

First, note that according to (\ref{B=Phi}) and (\ref{Lerch_def}),
\begin{equation}
\mathrm{B}\left( \nu ,0,z\right) =\sum_{k=0}^{\infty }\frac{z^{k+\nu }}{%
k+\nu },  \label{B_sum}
\end{equation}%
thus $\mathrm{B}\left( \nu ,0,z\right) $ is divergent for non-positive
integral values of $\nu $. Therefore, we will consider two different cases
in this Section:\ $\nu \in
%TCIMACRO{\U{211a} }%
%BeginExpansion
\mathbb{Q}
%EndExpansion
^{+}$ and $\nu \in
%TCIMACRO{\U{211a} }%
%BeginExpansion
\mathbb{Q}
%EndExpansion
^{-}\backslash \left\{ -1,-2,\ldots \right\} $.

\subsection{Case $\protect\nu \in
%TCIMACRO{\U{211a} }%
%BeginExpansion
\mathbb{Q}
%EndExpansion
^{+}$}

Consider $\nu =n+r>0$ where $n=\left\lfloor \nu \right\rfloor \geq 0$ is the
integer part of $\nu $ and $0\leq r\leq 1$. From (\ref{B_sum}), we have%
\begin{eqnarray}
\mathrm{B}\left( n+r,0,z\right)  &=&\sum_{k=0}^{\infty }\frac{z^{k+n+r}}{%
k+n+r}=\sum_{k=n}^{\infty }\frac{z^{k+r}}{k+r}  \notag \\
&=&\sum_{k=0}^{\infty }\frac{z^{k+r}}{k+r}-\sum_{k=0}^{n-1}\frac{z^{k+r}}{k+r%
}.  \label{Sums}
\end{eqnarray}

Set $r=1$ in (\ref{Sums}) and then apply the Taylor expansion \cite[Eqn.
4.6.1]{DLMF}
\begin{equation*}
\log \left( 1+z\right) =-\sum_{k=1}^{\infty }\frac{\left( -z\right) ^{k}}{k},
\end{equation*}%
to obtain%
\begin{eqnarray}
&&\mathrm{B}\left( n+1,0,z\right) =-\log \left( 1-z\right) -\sum_{k=1}^{n}%
\frac{z^{k}}{k},  \label{B(n+1)} \\
&&n=0,1,2,\ldots   \notag
\end{eqnarray}

Further, set $r=1/2$ in (\ref{Sums}) and apply de Taylor expansion \cite[%
Eqn. 4.38.5]{DLMF}
\begin{equation}
\tanh ^{-1}z=\sum_{k=0}^{\infty }\frac{z^{2k+1}}{2k+1},
\label{tanh-1_expansion}
\end{equation}%
to obtain%
\begin{eqnarray}
&&\mathrm{B}\left( n+\frac{1}{2},0,z\right) =2\left( \tanh ^{-1}\sqrt{z}%
-\sum_{k=0}^{n-1}\frac{z^{k+1/2}}{2k+1}\right) ,  \label{B(n+1/2)} \\
&&n=0,1,2,\ldots   \notag
\end{eqnarray}

More generally, set $r=p/q\in
%TCIMACRO{\U{211a} }%
%BeginExpansion
\mathbb{Q}
%EndExpansion
$ in (\ref{Sums}) with $p<q$. Then,
\begin{equation}
\mathrm{B}\left( n+\frac{p}{q},0,z\right) =z^{p/q}\sum_{k=0}^{\infty }\frac{%
z^{k}}{k+p/q}-\sum_{k=0}^{n-1}\frac{z^{k+p/q}}{k+p/q},  \label{B(n+p/q)_1}
\end{equation}

Rewrite the first sum of (\ref{B(n+p/q)_1})\ as a hypergeometric function
(see \cite[p. 61-62]{Andrews}),
\begin{equation}
\sum_{k=0}^{\infty }\frac{z^{k}}{k+p/q}=\frac{1}{p/q}\,_{2}F_{1}\left(
\left.
\begin{array}{c}
1,p/q \\
1+p/q%
\end{array}%
\right\vert z\right) .  \label{Sum=hyp}
\end{equation}

Apply now the reduction formula \cite[Eqn. 7.3.1.131]{Prudnikov3}%
\begin{eqnarray}
&&_{2}F_{1}\left( \left.
\begin{array}{c}
1,p/q \\
1+p/q%
\end{array}%
\right\vert z\right)  \notag \\
&=&-\frac{p}{q}z^{-p/q}\sum_{k=0}^{q-1}\exp \left( -\frac{2\pi ipk}{q}%
\right) \log \left( 1-z^{1/q}\exp \left( \frac{2\pi ik}{p}\right) \right) ,
\label{Hyp_p/p} \\
&&p,q=1,2,\ldots ;\ p\leq q.  \notag
\end{eqnarray}

Therefore, taking into account (\ref{Sum=hyp})-(\ref{Hyp_p/p}), rewrite (\ref%
{B(n+p/q)_1})\ as the following result.

\begin{theorem}
For $\nu =n+\frac{p}{q}\in
%TCIMACRO{\U{211a} }%
%BeginExpansion
\mathbb{Q}
%EndExpansion
^{+}$, with $n=\left\lfloor \nu \right\rfloor $ and $p<q$, the reduction
formula
\begin{eqnarray}
&&\mathrm{B}\left( \nu ,0,z\right) =z^{\nu }\mathrm{\Phi} \left( z,1,\nu \right)
\label{B_nu+} \\
&=&-\sum_{k=0}^{q-1}\exp \left( -\frac{2\pi ipk}{q}\right) \log \left(
1-z^{1/q}\exp \left( \frac{2\pi ik}{q}\right) \right) -\sum_{k=0}^{n-1}\frac{%
z^{k+p/q}}{k+p/q}.  \notag
\end{eqnarray}%
holds true.
\end{theorem}

\begin{remark}
Notice that the reduction formula (\ref{B(n+1/2)})\ is included in (\ref%
{B_nu+}), but not (\ref{B(n+1)}), which is a singular case.
\end{remark}

\subsection{Case $\protect\nu \in
%TCIMACRO{\U{211a} }%
%BeginExpansion
\mathbb{Q}
%EndExpansion
^{-}\backslash \left\{ -1,-2,\ldots \right\} $}

Consider $\nu =-n+r<0$ where $n=\left\lfloor \left\vert \nu -1\right\vert
\right\rfloor \geq 0$, and $0<r<1$. From (\ref{B_sum}), we have%
\begin{eqnarray}
\mathrm{B}\left( -n+r,0,z\right) &=&\sum_{k=0}^{\infty }\frac{z^{k-n+r}}{%
k-n+r}=\sum_{k=-n}^{\infty }\frac{z^{k+r}}{k+r}  \notag \\
&=&\sum_{k=0}^{\infty }\frac{z^{k+r}}{k+r}+\sum_{k=1}^{n}\frac{z^{-k+r}}{-k+r%
}.  \label{Sums_2}
\end{eqnarray}

Taking $r=1/2$ and applying again (\ref{tanh-1_expansion}), we have%
\begin{equation}
\mathrm{B}\left( -n+\frac{1}{2},0,z\right) =2\left( \tanh ^{-1}\sqrt{z}%
-\sum_{k=1}^{n}\frac{z^{-k+1/2}}{2k-1}\right) .  \label{B(-n+1/2)}
\end{equation}

More generally, take $r=p/q\in
%TCIMACRO{\U{211a} }%
%BeginExpansion
\mathbb{Q}
%EndExpansion
$ with $p<q$ in (\ref{Sums_2}), and apply (\ref{Sum=hyp}) to obtain the
following result.

\begin{theorem}
For $\nu =-n+\frac{p}{q}\in
%TCIMACRO{\U{211a} }%
%BeginExpansion
\mathbb{Q}
%EndExpansion
^{-}$, with $n=\left\lfloor \left\vert \nu -1\right\vert \right\rfloor $ and
$p<q$, the reduction formula
\begin{eqnarray}
&&\mathrm{B}\left( \nu ,0,z\right)   \label{B_nu-} \\
&=&-\sum_{k=0}^{q-1}\exp \left( -\frac{2\pi ipk}{q}\right) \log \left(
1-z^{1/q}\exp \left( \frac{2\pi ik}{q}\right) \right) +\sum_{k=1}^{n}\frac{%
z^{p/q-k}}{p/q-k}.  \notag
\end{eqnarray}%
holds true.
\end{theorem}

\begin{remark}
Notice that (\ref{B(-n+1/2)})\ is included in (\ref{B_nu-})\ as a particular
case. Also, in (\ref{B_nu-}), $\mathrm{B}\left( \nu ,0,z\right) \neq z^{\nu
}\mathrm{\Phi} \left( z,1,\nu \right) $ since (\ref{B=Phi})\ does not hold true for $%
\nu <0$.
\end{remark}

\section{Applications\label{Section: Applications}}

From the reduction formulas\ obtained in Section \ref{Section: Main Results}%
, next we calculate some integrals in terms of elementary functions. Also,
we will use these reduction formulas as a benchmark for the computation of
the incomplete beta function.

\subsection{Calculation of integrals}

On the one hand, an integral representation of the incomplete beta function
is given by \cite[Eqn. 58:3:5]{Atlas},
\begin{equation}
\mathrm{B}\left( \nu ,\mu ,z\right) =z^{\nu }\int_{0}^{1}t^{\nu -1}\left(
1-zt\right) ^{\mu -1}dt.  \label{beta_int}
\end{equation}

Also, an integral representation of the Lerch transcendent is \cite[Eqn.
1.11(3)]{Erdelyi}
\begin{equation}
\mathrm{\Phi} \left( z,s,\nu \right) =\frac{1}{\mathrm{\Gamma} \left( s\right) }%
\int_{0}^{\infty }\frac{t^{s-1}e^{-\left( \nu -1\right) t}}{e^{t}-z}dt,\quad
\mathrm{Re}\ \nu >0.  \label{Phi_int}
\end{equation}

Notice that from (\ref{B=Phi}), (\ref{beta_int}), and (\ref{Phi_int}), we
have%
\begin{eqnarray}
&&\int_{0}^{1}\frac{t^{\nu -1}}{1-zt}dt=\int_{0}^{\infty }\frac{e^{-\left(
\nu -1\right) t}}{e^{t}-z}dt=z^{-\nu }\mathrm{B}\left( \nu ,0,z\right) ,
\label{Int=} \\
&&\mathrm{Re}\ \nu >0.  \notag
\end{eqnarray}

Therefore, from (\ref{B_nu+})\ and (\ref{Int=}), we have for $\nu =n+\frac{p%
}{q}\in
%TCIMACRO{\U{211a} }%
%BeginExpansion
\mathbb{Q}
%EndExpansion
^{+}$, with $n=\left\lfloor \nu \right\rfloor $ and $p<q$,
\begin{eqnarray}
&&\int_{0}^{1}\frac{t^{\nu -1}}{1-zt}dt=\int_{0}^{\infty }\frac{e^{-\left(
\nu -1\right) t}}{e^{t}-z}dt  \label{New_ints} \\
&=&-z^{-n-p/q}\sum_{k=0}^{q-1}\exp \left( -\frac{2\pi ipk}{q}\right) \log
\left( 1-z^{1/q}\exp \left( \frac{2\pi ik}{q}\right) \right)
-\sum_{k=0}^{n-1}\frac{z^{k-n}}{k+p/q}.  \notag
\end{eqnarray}

On the other hand, in the literature we found \cite[Eqn. 58:14:7]{Atlas}%
\begin{equation}
\int_{0}^{z}\tanh ^{2\lambda -1}t\ dt=\frac{1}{2}\mathrm{B}\left( \lambda
,0,\tanh ^{2}z\right) ,\quad \mathrm{Re}\ \lambda >0.  \label{Int_tanh}
\end{equation}

Therefore, from (\ref{B_nu+})\ and (\ref{Int_tanh}), we have for $\lambda =n+%
\frac{p}{q}\in
%TCIMACRO{\U{211a} }%
%BeginExpansion
\mathbb{Q}
%EndExpansion
^{+}$, with $n=\left\lfloor \lambda \right\rfloor $ and $p<q$,%
\begin{eqnarray}
&&\int_{0}^{z}\tanh ^{2\lambda -1}t\ dt  \label{Int_tanh_resultado} \\
&=&-\frac{1}{2}\sum_{k=0}^{q-1}\exp \left( -\frac{2\pi ipk}{q}\right) \log
\left( 1-\left( \tanh z\right) ^{2/q}\exp \left( \frac{2\pi ik}{q}\right)
\right)   \notag \\
&&-\frac{1}{2}\sum_{k=0}^{n-1}\frac{\left( \tanh z\right) ^{2\left(
k+p/q\right) }}{k+p/q}.  \notag
\end{eqnarray}

The integral given in (\ref{Int_tanh_resultado})\ generalizes the results
found in the literature for $\lambda =n+1$ and $\lambda =n+\frac{1}{2}$ with
$n=0,1,2,\ldots $ \cite[Eqns. 2.424.2-3]{Gradstheyn}.

It is worth noting that for particular values of $\lambda $, the \textsf{%
Integrate} MATHEMATICA command is able to compute symbolically
the same results as (\ref{Int_tanh_resultado}), but in a very time-consuming
way. For instance, for $\lambda =\frac{5}{4}$, we obtain%
\begin{equation*}
\int_{0}^{z}\tanh ^{3/2}t\ dt=\tanh ^{-1}\left( \sqrt{\tanh z}\right) -2%
\sqrt{\tanh z}+\tan ^{-1}\left( \sqrt{\tanh z}\right) ,
\end{equation*}%
but the \textsf{Integrate} command takes around $300$ times longer than the
reduction formula given in (\ref{Int_tanh_resultado}).

\subsection{Numerical evaluation}

From a numerical point of view, the reduction formulas (\ref{B_nu+})\ and (%
\ref{B_nu-}) are quite useful to plot $\mathrm{B}\left( \nu ,0,z\right) $ as
a function of $\nu $ in the real domain. However, for some real values of $%
\nu $ and $z$, we obtain a complex value for $\mathrm{B}\left( \nu
,0,z\right) $. In these cases, the imaginary part of $\mathrm{B}\left( \nu
,0,z\right) $ is not always easy to compute. Figure \ref{Figure: Numerical
difference} shows the plot of $\mathrm{Im}\left( \mathrm{B}\left( \nu
,0,z\right) \right) $ as a function of $z$ for $\nu =12.3$. The reduction
formula (\ref{B_nu+}) shows the correct answer, i.e. $\mathrm{Im}\left(
\mathrm{B}\left( \nu ,0,z\right) \right) =-\pi $, meanwhile the numerical
evaluation of $\mathrm{Im}\left( \mathrm{B}\left( \nu ,0,z\right) \right) $
with MATHEMATICA\ diverges from this result. A\ similar
feature is observed using (\ref{B_nu-}) and a negative value for $\nu $. It
is worth noting that the equivalent numerical evaluation of $\mathrm{Im}%
\left( z^{\nu }\mathrm{\Phi} \left( z,1,\nu \right) \right) $ with MATHEMATICA
yields also $-\pi $. \ 

\begin{figure}
  \includegraphics{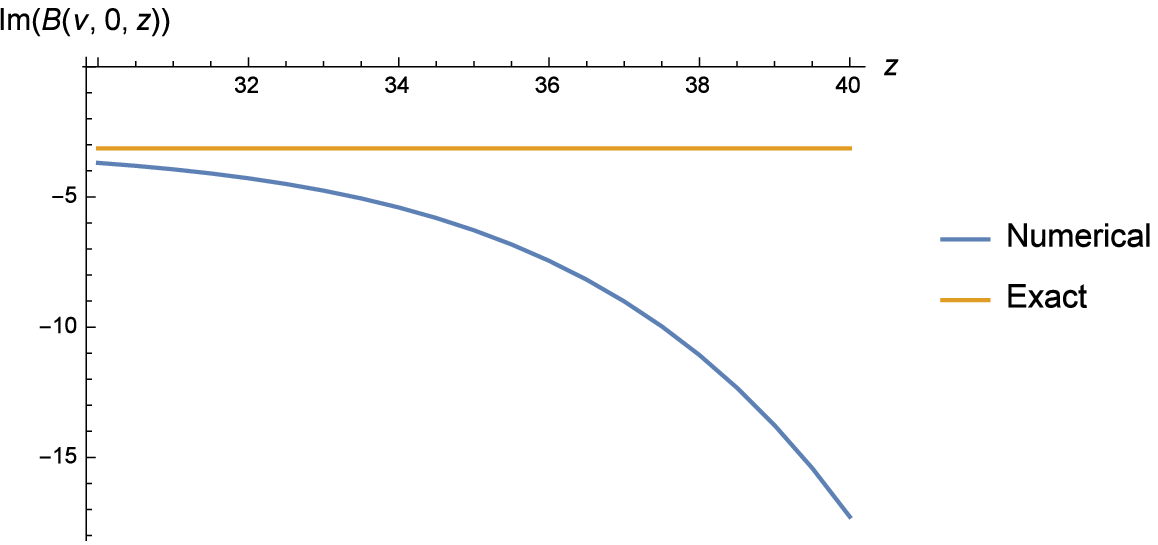}
\caption{Evaluation of $\mathrm{Im}\left( \mathrm{B}\left( \protect\nu %
,0,z\right) \right) $ with MATHEMATICA\ and (\protect\ref%
{B_nu+}) with $\protect\nu =12.3$.}
\label{Figure: Numerical difference}       % Give a unique label
\end{figure}
\section{Conclusions\label{Section: Conclusions}}

On the one hand, we have derived in (\ref{B_nu+})\ and (\ref{B_nu-})\ new
expressions for the incomplete beta function $\mathrm{B}\left( \nu
,0,z\right) $ and the Lerch transcendent $\mathrm{\Phi} \left( z,1,\nu \right) $\ in
terms of elementary functions when $\nu $ is rational and $z$ is complex.
Particular formulas for non-negative integers values of $\nu $ and for
half-integer values of $\nu $ are given in (\ref{B(n+1)})\ and (\ref%
{B(n+1/2)}), (\ref{B(-n+1/2)}) respectively.

On the other hand, we have calculated the integrals given (\ref{New_ints})\
from the reduction formulas (\ref{B_nu+})\ and (\ref{B_nu-}) and the
integral representation of the incomplete beta function and the Lerch
transcendent. Also, in (\ref{Int_tanh_resultado}), the integral $%
\int_{0}^{z}\tanh ^{\alpha }t\ dt$ is calculated in terms of elementary
functions for $\alpha \in
%TCIMACRO{\U{211a} }%
%BeginExpansion
\mathbb{Q}
%EndExpansion
$ and $\alpha >-1$. It is worth noting that (\ref{Int_tanh_resultado})
accelerates quite significantly the symbolic computation of the latter
integral with the aid of computer algebra.

Finally, with the aid of the reduction formulas (\ref{B_nu+})\ and (\ref%
{B_nu-}), we have tested that the numerical algorithm provided by MATHEMATICA 
sometimes fails to compute the imaginary part of $\mathrm{B}%
\left( \nu ,0,z\right) $. Also, the reduction formulas (\ref{B_nu+})\ and (%
\ref{B_nu-})\ are numerically useful to plot $\mathrm{B}\left( \nu
,0,z\right) $ as a function of $\nu $ in the real domain.

All the results presented in this paper have been implemented in MATHEMATICA 
and can be downloaded from https://bit.ly/2XT7UjK

%\begin{acknowledgements}
%\end{acknowledgements}

\section*{Conflict of interest}
The authors declare that they have no conflict of interest.



\end{document}